\documentclass[11pt]{article}
\usepackage{amsmath,amssymb,eucal}
\usepackage[dvips]{graphicx}

\newtheorem{thm}{\bf Theorem}[section]

\numberwithin{equation}{section}

\def\qed{$\Box$}

 \newcommand{\mapdown}[1]{\Big\downarrow
\llap{$\vcenter{\hbox{$\scriptstyle#1\,$}}$}}
 
\begin{document}

\title{{\bf 
{\large Seifert surgery on knots via Reidemeister torsion and Casson-Walker-Lescop invariant III}}}
\author{
{\normalsize Teruhisa Kadokami, Noriko Maruyama and Tsuyoshi Sakai}}
\date{{\normalsize March 16th, 2018}}
\footnotetext[0]{%
2010 {\it Mathematics Subject Classification}:
11R04, 11R27, 57M25, 57M27. \par
{\it Keywords}: 
Reidemeister torsion,
Casson-Walker-Lescop invariant,
Seifert fibered space.}
\maketitle

\begin{abstract}{
For a knot $K$ in a homology $3$-sphere $\Sigma$,
let $M$ be the result of $2/q$-surgery on $K$, and
let $X$ be the universal abelian covering of $M$.
Our first theorem is that if the first homology of $X$ is finite cyclic
and $M$ is a Seifert fibered space with $N\ge 3$ singular fibers, 
then $N\ge 4$ if and only if 
the first homology of the universal abelian covering of $X$ is infinite.
Our second theorem is that under an appropriate assumption 
on the Alexander polynomial of $K$, if $M$ is a Seifert fibered space, 
then $q=\pm 1$ (i.e.\ integral surgery).
}\end{abstract}

\section{Introduction}\label{sec:intro}

It is conjectured that Seifert surgeries on non-trivial knots are integral (except some cases), 
and that the resultant Seifert fibered spaces have three singular fibers
(for example, see \cite{DMM}).
We \cite{KMS1, KMS2} have studied the $2/q$-Seifert surgeries.
One reason why the coefficients are $2/q$ is that we started from the following 
shown in \cite{Kd1} : 

\begin{thm}\label{th:Seifert}
{\rm (\cite[Theorem 1.4]{Kd1})}
Let $K$ be a knot in a homology $3$-sphere $\Sigma$ such that
the Alexander polynomial of $K$ is $t^2-3t+1$.
The only surgeries on $K$ that may produce a Seifert fibered space with base $S^2$
and with $H_1\ne \{0\}, \mathbb{Z}$ have coefficients $2/q$ and $3/q$,
and produce Seifert fibered space with three singular fibers.
Moreover
(1) if the coefficient is $2/q$, then the set of multiplicities is
$\{2\alpha, 2\beta, 5\}$ where $\gcd(\alpha, \beta)=\gcd(\alpha, 5)=\gcd(\beta, 5)=1$, and
(2) if the coefficient is $3/q$, then the set of multiplicities is
$\{3\alpha, 3\beta, 4\}$ where $\gcd(\alpha, \beta)=\gcd(\alpha, 2)=\gcd(\beta, 2)=1$.
\end{thm}

In this paper, we continue to study the $2/q$-Seifert surgeries, 
where we replace the condition on the Alexander polynomial 
with a condition on the first homology of the double branched covering.
We note that the results below are related to the conjectures 
mentioned above.

\medskip

\noindent
(1) Let $\Sigma$ be a homology $3$-sphere, 
and let $K$ be a knot in $\Sigma$.
Then 
$\Delta_K(t)$ denotes the Alexander polynomial of $K$, 
$a_2(K)$ denotes the Conway's $a_2$-term of $K$, 
and $\Sigma(K; p/r)$ denotes the result of $p/r$-surgery on $K$.

\bigskip

\noindent
(2) The first author \cite{Kd2} introduced the norm of polynomials and homology lens spaces: 
Let $\zeta_d$ be a primitive $d$-th root of unity.
For an element $\alpha$ of $\mathbb{Q}(\zeta_d)$,
$N_d(\alpha)$ denotes the norm of $\alpha$ associated to
the algebraic extension $\mathbb{Q}(\zeta_d)$ over $\mathbb{Q}$.
Let $f(t)$ be a Laurent polynomial over $\mathbb{Z}$.
We define $|f(t)|_d$ by
$$|f(t)|_d=|N_d(f(\zeta_d))|
=\left| \prod_{i\in (\mathbb{Z}/d\mathbb{Z})^{\times}}
f(\zeta_d^i)\right|.$$
Let $X$ be a homology lens space with $H_1(X)\cong \mathbb{Z}/p\mathbb{Z}$.
Then there exists a knot $K$ in a homology $3$-sphere $\Sigma$ such that
$X=\Sigma(K; p/r)$ (\cite[Lemma 2.1]{BL}).
We define $|X|_d$ and $\|X\|_d$ by
$$|X|_d=|\Delta_K(t)|_d\quad
\mbox{and}\quad
\|X\|_d=\prod_{d'|d, d'\ge 1}|X|_{d'},$$
where $d$ is a divisor of $p$.
Then both $|X|_d$ and $\|X\|_d$ are topological invariants of $X$ (Refer to \cite{Kd2} for details).

\bigskip

\noindent
(3) Let $X$ be a closed oriented $3$-manifold.
Then $\lambda(X)$ denotes the Lescop invariant of $X$ (\cite{Ls}). 
Note that $\lambda\left(S^3\right)=0$.

\section{Results}\label{sec:result}

Let $K$ be a knot in a homology $3$-sphere $\Sigma$.
For an odd integer $q$, 
let $M$ be the result of $2/q$-surgery on $K$: $M=\Sigma(K; 2/q)$.
Let $\pi : X\to M$ be the universal abelian covering of $M$
(i.e.\ the covering associated to $\mathrm{Ker}(\pi_1(M)\to H_1(M))$).
Since $H_1(M)\cong \mathbb{Z}/2\mathbb{Z}$, $\pi$ is the $2$-fold unbranched covering.

\medskip

Let $\Sigma_2$ be the double branched covering space of $\Sigma$
branched along $K$, and $\overline{K}$ the lifted knot of $K$ in $\Sigma_2$.
Since $\overline{K}$ is null-homologous in $\Sigma_2$, and $X$ is 
the result of $1/q$-surgery on $\overline{K}$, we have 
$H_1(X)\cong H_1(\Sigma_2)$.

$$
\begin{array}{rcccl}
\overline{K}\subset & \Sigma_2 & \leadsto & X & =\Sigma_2(\overline{K}; 1/q)\\
& \mapdown{} & & \mapdown{} & \\
K\subset & \Sigma & \leadsto & M & =\Sigma(K; 2/q)
\end{array}$$

\bigskip

Assume that 
$H_1(\Sigma_2)\cong \mathbb{Z}/m\mathbb{Z}$ with $m\ge 3$.
Then $H_1(X)\cong \mathbb{Z}/m\mathbb{Z}$.
Hence $|X|_d$ is defined for each divisor $d$ of $m$, and 
so is $\|X\|_d$.

\medskip

We then have the following.

\begin{thm}\label{th:main1}
Let $K$ be a knot in a homology 3-sphere $\Sigma$. 
With the notation above, 
we assume that $H_1(\Sigma_2)\cong \mathbb{Z}/m\mathbb{Z}$ with $m\ge 3$, 
and $M=\Sigma(K; 2/q)$ is a Seifert fibered space with $N$ singular fibers 
where $N\ge 3$.
Then the following hold.

\medskip
\noindent
(1) $N\ge 4$ if and only if $\|X\|_m=0$.

\medskip
\noindent
(2) If $m$ is a power of a prime, 
then $N=3$ and the set of multiplicities is $\{2\alpha, 2\beta, m\}$, 
where $\alpha$, $\beta$ and $m$ are pairwise coprime.

\end{thm}

\newpage

\begin{thm}\label{th:main2}
Let $K$ be a knot in a homology 3-sphere $\Sigma$. We assume the following.
 
\medskip
\noindent
{\bf (2.1)} $\lambda(\Sigma)=0$.

\medskip
\noindent
{\bf (2.2)} $|\Delta_K(-1)|=5$.
 
\medskip
\noindent
{\bf (2.3)} $|q|\ge 3$, 
$\sqrt{|X|_5}>4\{qa_2(K)\}^2$.

\medskip

Then $M=\Sigma(K; 2/q)$ is not a Seifert fibered space.

\end{thm}

\section{Proof of Theorem \ref{th:main1}}\label{sec:proof1}

Since $H_1(M)\cong \mathbb{Z}/2\mathbb{Z}$, 
the base surface of $M$ is $S^2$.
As noted in Section 2, 
$H_1(X)\cong H_1(\Sigma_2)\cong \mathbb{Z}/m\mathbb{Z}$.
Hence $\|X\|_m$ is defined, and also we have
$$\|M\|_2=|M|_2=|\Delta_K(-1)|=|H_1(\Sigma_2)|=m\ge 3.$$
Hence, by \cite[Theorem 1.3 (2)]{Kd1}, 
we may assume that $M$ has a framed link presentation as in Figure 1, 
where $\gcd(p_1, p_2)\ge 2$,  
$p_k\ge 2$ for $k=1, \ldots, N$, 
and $\gcd(p_i, p_j)=1$ for $1\le i< j\le N$ and $\{i, j\}\ne \{1, 2\}$.

\begin{figure}[htbp]
\begin{center}
\includegraphics[scale=0.5]{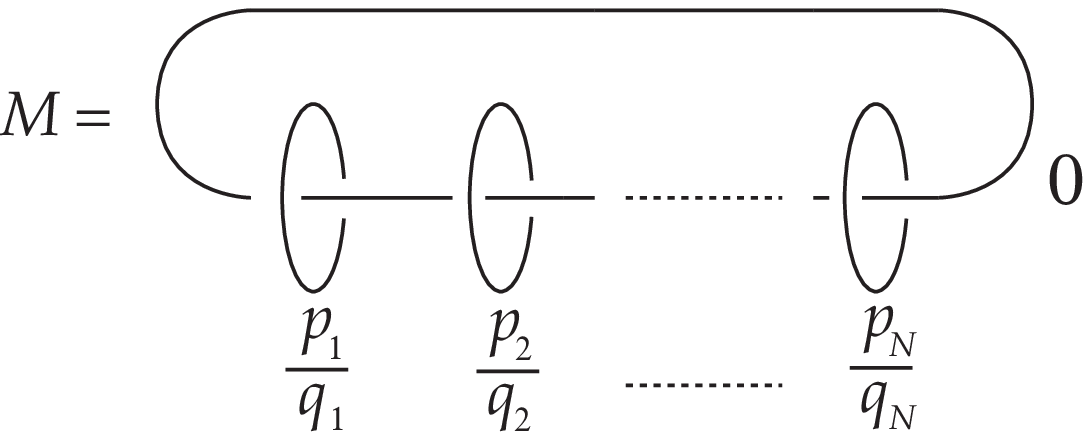}
\label{M1}
\caption{A framed link presentation of $M=\Sigma (K; 2/q)$}
\end{center}
\end{figure}

Since $H_1(M)\cong \mathbb{Z}/2\mathbb{Z}$, 
we have $\gcd(p_1, p_2)= 2$ by \cite[Theorem 1.2 (1)]{Kd1}.
Hence we may assume that 
$M$ has a framed link presentation as in Figure 2, 
where $\alpha$, $\beta$, $p_3, \ldots, p_N$ are pairwise coprime 
and $p_3, \ldots, p_N$ are all odd.

\begin{figure}[htbp]
\begin{center}
\includegraphics[scale=0.5]{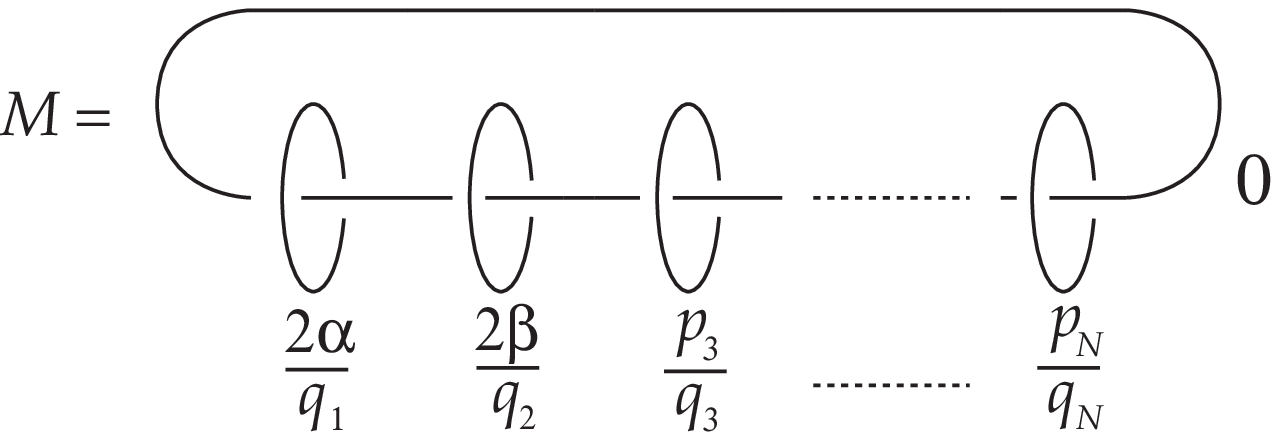}
\label{M2}
\caption{A framed link presentation of $M=\Sigma (K; 2/q)$}
\end{center}
\end{figure}

Hence, as in \cite{KMS1}, we see that 
$X$ has a framed link presentation as in Figure 3.

\newpage

\begin{figure}[htbp]
\begin{center}
\includegraphics[scale=0.5]{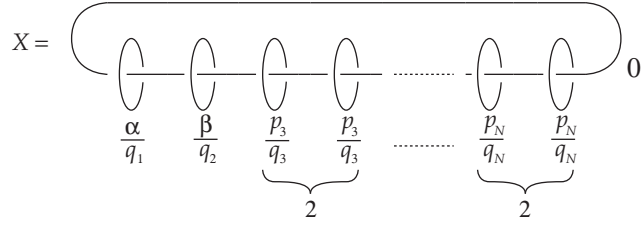}
\label{X}
\caption{A framed link presentation of $X$}
\end{center}
\end{figure}

Hence, by \cite[Theorem 1.3]{Kd1}, we have (1).

\medskip

Moreover, by \cite[Theorem 1.2 (3)]{Kd1}, we have 
$|M|_2=p_3\cdots p_N$, and hence $p_3\cdots p_N=m$.
Therefore if $m$ is a power of a prime, we have $N=3$ and $p_3=m$.
This completes the proof.
\qed


\section{Proof of Theorem \ref{th:main2}}\label{sec:proof2}

Suppose that $M=\Sigma(K; 2/q)$ is a Seifert fibered space.
Recall that 
$|H_1(\Sigma_2)|=|\Delta_K(-1)|=5$.
Then, by Theorem 2.1 (2), we may assume that
$M$ has a framed link presentation as in Figure 4, where 
$1\le \alpha<\beta$ and $\gcd(\alpha, \beta)=1$.

\begin{figure}[htbp]
\begin{center}
\includegraphics[scale=0.5]{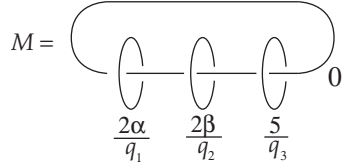}
\label{M2}
\caption{A framed link presentation of $M=\Sigma (K; 2/q)$}
\end{center}
\end{figure}

\medskip

Since $\Delta_K(-1)\equiv \pm 3\ (\mathrm{mod}\ \! 8)$ 
according to {\bf (2.2)}, 
\cite[Theorem 10.7, p108]{Lc} implies that $a_2(K)\ne 0$.

\medskip

As shown in \cite[Theorem 1.2 (3)]{Kd1}, $\sqrt{|X|_5}=(\alpha \beta)^2$.
By {\bf (2.1)} and \cite[1.5 T2]{Ls}, we have $\lambda(M)=qa_2(K)$.
By substituting in $\sqrt{|X|_5}>4\{qa_2(K)\}^2$, 
we have $(\alpha \beta)^2>4\{\lambda(M)\}^2$.
Hence we have 
$|\lambda(M)|<(\alpha \beta)/2$.
The remaining part of the proof goes parallel to the proof of 
\cite[Theorem 2.1]{KMS2}.
This completes the proof.
\qed

\bigskip

\noindent
{\bf Acknowledgement}\ 
The first author is supported by Grant-in-Aid for Scientific Research(C) 
Number 17K05246, and 
the second author is supported by Grant-in-Aid for Scientific Research(C) 
Number 16K05158.

{\small
\par
Teruhisa Kadokami\par
School of Mechanical Engineering,\par
College of Science and Engineering, Kanazawa University,\par
Kakuma-machi, Kanazawa, Ishikawa, 920-1192, Japan\par
{\tt kadokami@se.kanazawa-u.ac.jp}\par

\medskip

Noriko Maruyama\par
Musashino Art University,\par 
Ogawa 1-736, Kodaira, Tokyo 187-8505, Japan \par 
{\tt maruyama@musabi.ac.jp} \par

\medskip

Tsuyoshi Sakai\par
Department of Mathematics, Nihon University,\par
3-25-40, Sakurajosui, Setagaya-ku, Tokyo 156-8550, Japan \par
}


\newpage

{\footnotesize
 \begin{thebibliography}{999}
%
\newcommand{\bysame}{%
       \leavevmode\hbox to 3em{\hrulefill}\,}
%
\bibitem[BL]{BL}S.~Boyer and D.~Lines,
{\it Surgery formulae for Casson's invariant and extensions to homology lens spaces},
J. Reine Angew. Math., {\bf 45}\ (1990), 181--220.
%
\bibitem[DMM]{DMM}A.~Deruelle, K.~Miyazaki and K.~Motegi,
{\it Networking Seifert surgeries on knots},
Mem. Amer. Math. Soc. {\bf 217} (2012), no. 1021.
%
\bibitem[Kd1]{Kd1}T.~Kadokami,
{\it Reidemeister torsion of Seifert fibered homology lens spaces and Dehn surgery},
Algebr. Geom. Topol., {\bf 7}\ (2007), 1509--1529.
%
\bibitem[Kd2]{Kd2}T.~Kadokami,
{\it Reidemeister torsion and lens surgeries on knots in homology $3$-spheres II},
Top. Appl., {\bf 155}, no.15\ (2008), 1699--1707.
%
\bibitem[KMS1]{KMS1}T.~Kadokami, N.~Maruyama and T.~Sakai, 
{\it Seifert surgery on knots via Reidemeister torsion and Casson-Walker-Lescop invariant}, 
Top. Appl., {\bf 188} (2015), 64--73.
%
\bibitem[KMS2]{KMS2}T.~Kadokami, N.~Maruyama and T.~Sakai, 
{\it Seifert surgery on knots via Reidemeister torsion and Casson-Walker-Lescop invariant II}, 
Osaka J. Math., {\bf 53} (2016), 767--773.
%
\bibitem[Ls]{Ls}C.~Lescop, 
{\it Global surgery formula for the Casson-Walker invariant},
Ann. of Math. Studies, Princeton Univ. Press., {\bf 140} (1996).
%
%
\bibitem[Lc]{Lc}W.~B.~R.~Lickorish, 
{\it An introduction to knot theory},
Graduate texts in mathematics, {\bf 175}\ (1997).
%
\end{thebibliography}}
\end{document}